\documentclass[a4paper,12pt]{article}
\usepackage[cp1251]{inputenc}
\usepackage[russian, ukrainian,english]{babel}
\usepackage{amsmath, amsthm, amsfonts, amssymb}

\theoremstyle{plain}
\newtheorem{theorem}{Theorem}[section]

\newtheorem{lemma}{Lemma}[section]

\begin{document}

\begin{center}
{\bf\Large Krylov-Veretennikov Formula for Functionals from the Stopped Wiener Process}
\vskip 10 pt
{\bf G. V. Riabov}
\vskip 10 pt
{\it Institute of Mathematics, NAS of Ukraine}
\end{center}

\vskip 10 pt
{\bf Abstract}
\vskip 10 pt
We consider a class of measures absolutely continuous with respect to the distribution of the stopped Wiener process $w(\cdot\wedge\tau)$. Multiple stochastic integrals, that lead to the analogue of the It\^o-Wiener expansions for such measures, are described. An analogue of the Krylov-Veretennikov formula for functionals $f=\varphi(w(\tau))$ is obtained.

\section{Introduction}

Let $\{w(t)\}_{t\geq 0}$ be a standard Wiener process in $\mathbb{R}^d,$ starting from the point $u\in \mathbb{R}^d.$ Consider an open connected set  $G\ni u,$ the exit time
$$
\tau=\inf\{t>0: w(t)\not\in G\},
$$
and a Borel function $\rho:\mathbb{R}^d\to (0,1).$

The main object of the investigation in the present paper is the orthogonal structure of the space $L^2(\Omega,\sigma(w(\cdot\wedge \tau)),Q),$ where the measure $Q$ is given by the  density
$$
\frac{dQ}{d\mathbb{P}}=\frac{1_{\tau<\infty}\rho(w(\tau))}{\mathbb{E}1_{\tau<\infty}\rho(w(\tau))}.
$$

In \cite[L. 2.4]{Riabov} it was proved that the space $L^2(\Omega,\sigma(w(\cdot\wedge \tau)),Q)$ possesses an orthogonal structure similar to the It\^o-Wiener decomposition in the Gaussian case \cite{CM, Ito, Nualart}. Namely, consider functions
$$
\beta(v)=\mathbb{E} 1_{\tau(w-u+v)<\infty}\rho(w(\tau(w-u+v)-u+v)), \ v\in G,
$$
$$
\alpha(s,v)=\beta^{-1}(v)\mathbb{E} 1_{s<\tau(w-u+v)<\infty}\rho(w(\tau(w-u+v)-u+v)), \ s>0, v\in G,
$$
and processes

$$
\tilde{w}(s)=w(s\wedge\tau)-\int^{s\wedge\tau}_0\nabla \log\beta(w(r))dr, \ s\geq 0;
$$
$$
\hat{\tilde{w}}_t(s)=\tilde{w}(s)-\int^{s\wedge\tau}_0\nabla \log\alpha(t-r,w(r))dr, \ 0\leq s\leq t.
$$

\begin{theorem}
	\label{thm_IW}
	\cite[L. 2.4]{Riabov} Each random variable $f\in L^2(\Omega,\sigma(w(\cdot\wedge\tau)),Q)$ can be uniquely represented as a series of pairwise orthogonal stochastic integrals
	\begin{equation}
		\label{intro_IW_stopped_density}
		f=\sum^\infty_{n=0}\underset{0<t_1<\ldots<t_n<\tau}{\int\ldots\int}a_n(t_1,\ldots,t_n)d\hat{\tilde{w}}_{t_n}(t_1)\ldots d\hat{\tilde{w}}_{t_n}(t_{n-1})d\tilde{w}(t_n).
	\end{equation}
	Conversely, given a sequence of Borel functions $a_n:(0,\infty)^n\to \mathbb{R}^{d^n},$ $n\geq 0,$ such that
	$$
	\sum^\infty_{n=0}\underset{0<t_1<\ldots<t_n}{\int\ldots\int}\alpha(t_n,u) |a_n(t_1,\ldots,t_n)|^2dt_1\ldots dt_n<\infty,
	$$
	the series in the right-hand side of \eqref{intro_IW_stopped_density} converges in $L^2(\Omega,\sigma(w(\cdot\wedge\tau)),Q),$ and its sum $f$ satisfies
	$$
	\mathbb{E}f^2=\sum^\infty_{n=0}\underset{0<t_1<\ldots<t_n}{\int\ldots\int}
	\alpha(t_n,u) |a_n(t_1,\ldots,t_n)|^2dt_1\ldots dt_n.
	$$
\end{theorem}

In this paper we derive the explicit form of the expansion \eqref{intro_IW_stopped_density} for random variables of the kind $f=\varphi(w(\tau)).$ The resulting formula is similar to the well-known Krylov-Veretennikov formula \cite{KV}. It is written in terms of the transition semigroup $\{T^k_t\}_{t\geq 0}$ of a certain diffusion, killed at the boundary of $G.$ Indeed, the process $\tilde{w}$ is a stopped Wiener process relatively to the measure $Q$ \cite[L. 2.4]{Riabov}. Respectively, the initial Wiener process $w$ is a diffusion process relatively to the measure $Q.$ Then $\{T^k_t\}_{t\geq 0}$ is the transition semigroup of the process $w$ killed at the boundary of $G.$ Let $T$ denote the integration with respect to the exit distribution of $w$ from $G$ (precise expressions for these operators are given in the section 2). The main result of the present paper is the following formula, proved in the theorem \ref{KV_stopped}: 

	for every random variable $\varphi(w(\tau))\in L^2(\Omega,\sigma(w(\cdot\wedge\tau)),Q)$ the expansion \eqref{intro_IW_stopped_density} has the form
	\begin{equation}
	\label{intro_KV_stopped_eq}
	\begin{gathered}
	\varphi(w(\tau))=\sum^\infty_{n=0}
	\underset{0<t_1<\ldots<t_n}{\int\ldots\int}\alpha(t_n,u)^{-1} \bigg(T^k_{t_1}\alpha(t_2-t_1,\cdot)\nabla\big(\alpha(t_2-t_1,\cdot)^{-1} T^k_{t_2-t_1}\big)\ldots \\
	\alpha(t_n-t_{n-1},\cdot)\nabla\big(\alpha(t_n-t_{n-1},\cdot)^{-1}T^k_{t_n-t_{n-1}}\big) \nabla T\varphi\bigg)(u)
	d\hat{\tilde{w}}_{t_n}(t_1)\ldots d\hat{\tilde{w}}_{t_n}(t_{n-1})d\tilde{w}(t_{n}).
	\end{gathered}
	\end{equation}

Expansions of the kind \eqref{intro_IW_stopped_density} appeared in \cite{DorRiabov} in connection with the problem of studying the behaviour of Gaussian measures under nonlinear transformations. Such expansions have two main features:

\begin{enumerate}
	\item the summands in \eqref{intro_IW_stopped_density} are pairwise orthogonal;
	
	\item the summands in \eqref{intro_IW_stopped_density} are $\sigma(w(\cdot\wedge\tau))-$measurable.
\end{enumerate}
Of course, there are other possibilities to organize series expansions for random variables from $L^2(\Omega,\sigma(w(\cdot\wedge \tau)),Q).$ For simplicity, consider the case $Q=\mathbb{P}.$ The most straightforward approach comes from the obvious inclusion $\sigma(w(\cdot\wedge\tau))\subset \sigma(w).$ It means that each random variable
$f\in L^2(\Omega,\sigma(w(\cdot\wedge \tau)),\mathbb{P})$ possesses an It\^o-Wiener expansion with respect to the Wiener process $w:$
\begin{equation}
\label{IW}
f=\sum^\infty_{n=0}\underset{0<t_1<\ldots<t_n}{\int\ldots\int}b_n(t_1,\ldots,t_n)dw(t_1)\ldots dw(t_n).
\end{equation}
The summands in the expansion are not $\sigma(w(\cdot\wedge\tau))-$measurable. While the left-hand side of \eqref{IW} is $\sigma(w(\cdot\wedge\tau))-$measurable, one can condition \eqref{IW} with respect to $w(\cdot\wedge\tau)$ and get another expansion
\begin{equation}
\label{IW_conditioned}
f=\sum^\infty_{n=0}\underset{0<t_1<\ldots<t_n<\tau}{\int\ldots\int}b_n(t_1,\ldots,t_n)dw(t_1)\ldots dw(t_n).
\end{equation}
Now the stochastic integrals of different degree are not orthogonal. This causes known inconveniences: the expansion \eqref{IW_conditioned} is not unique (an example is given in \cite{DorRiabov}); the conditions for the expression in the right-hand side of \eqref{IW_conditioned} to converge are complicated. An application of the Gram-Shmidt orthogonalization procedure to expansions \eqref{IW_conditioned} was considered in \cite{DorIW}. However, in our framework it seems to be too complicated either to obtain the orthogonalized form of \eqref{IW_conditioned}, or to find the orthogonalized expansion \eqref{IW_conditioned} for a concrete random variable $f.$ The expansion \eqref{intro_IW_stopped_density} overcomes all these problems.

Motivation for the $\sigma(w(\cdot\wedge\tau))-$measurability of the summands in \eqref{intro_IW_stopped_density} comes from B. S. Tsirelson's theory of black noise. It is well-known that Brownian coalescing flows produce filtrations with trivial Gaussian parts \cite{Tsirelson, LJR}. So, to get a unified description of functionals measurable with respect to such flows, it is reasonable to use the noise generated by the flow itself. The results from \cite{DorRiabov, Riabov} show that this idea works: in \cite{DorRiabov} an orthogonal expansion of the kind \eqref{intro_IW_stopped_density} was obtained for the stopped Brownian motion; in \cite{Riabov} the same was done for the $n-$point motions of the Arratia flow. We refer to \cite{DorRiabov, Riabov} for the detailed discussion of this and related questions.

Generalization of the Krylov-Veretennikov formula to the wide class of dynamical systems driven by the additive Gaussian noise was obtained in \cite{Dor_KV}. Our formula \eqref{KV_stopped_eq} is similar to the one obtained in \cite{Dor_KV} despite the additional multipliers $\alpha.$ They occure to normalize operators $T^k_t,$ as $T^k_t1=\alpha(t,\cdot).$ 

The article is organized in the following way. In the section 2 we introduce all the needed notions and constructions. Also, it contains the reduction of the main theorem \ref{KV_stopped} to lemmata \ref{Clark_stopped} and \ref{KV_killed}. Sections 3 and 4 are devoted to the proof of these auxiliary results.

\section{Notations and Main Results}

To formulate our results, we will use the following notations.

$\{w(t)\}_{t\geq 0}$ is the Wiener process in $\mathbb{R}^d.$ Without loss of generality, we will assume that $w$ is constructed in a canonical way:

$\Omega=C([0,\infty),\mathbb{R}^d)$ is a space of continuous functions equipped with a metric of uniform convergence on compacts;

$\mathcal{F}$ is the Borel $\sigma-$field on $\Omega;$

$w(t,\omega)=\omega(t)$ is the canonical process on $(\Omega,\mathcal{F}),$  $\mathcal{F}_t=\sigma(w(s):0\leq s\leq t)$ is the natural filtration of $w;$

$(\mathbb{P}_v)_{v\in \mathbb{R}^d}$ is a family of probability measures on $(\Omega,\mathcal{F}),$ such that relatively to $\mathbb{P}_v,$  $w$ is a $d-$dimensional Wiener process starting from $v.$  The expectation with respect to certain probability measure $Q$ on $(\Omega,\mathcal{F})$ will be denoted by $\mathbb{E}_Q.$ $\mathbb{E}_{\mathbb{P}_v}$ will be abbreviated to $\mathbb{E}_v.$

Let $G\subset \mathbb{R}^d$ be an open connected set, $\tau$  be the exit time of $w$ from the set $G:$
$$
\tau=\inf\{t>0: w(t)\not \in G\}.
$$
We will assume that for all $v\in G,$  $\mathbb{P}_v(\tau<\infty)>0.$ Fix a Borel function $\rho:\mathbb{R}^d\to (0,1)$ and consider the function
$$
\beta(v)=\mathbb{E}_v 1_{\tau<\infty}\rho(w(\tau)), \ v\in G.
$$
It is a harmonic function in $G$ \cite[Ch. 4, Prop. 2.1]{Port_Stone}:
$$
\Delta_v \beta(v) =0, \ v\in G.
$$

Denote $Q_u$ the probability measure on $(\Omega,\mathcal{F}),$ defined via the density
$$
\frac{dQ_u}{d\mathbb{P}_u}=\beta(u)^{-1}1_{\tau<\infty}\rho(w(\tau)).
$$

We will need another probability measure corresponding to the process $w$ killed at the moment $\tau.$ Consider the function
$$
\alpha(s,v)=Q_v(\tau>s), \ s>0, v\in G.
$$
In the section 1 following processes were introduced.
\begin{equation}
\label{tilde_w}
\tilde{w}(s)=w(s\wedge\tau)-\int^{s\wedge\tau}_0\nabla_v \log\beta(w(r))dr, \ s\geq 0;
\end{equation}

\begin{equation}
\label{hat_tilde_w}
\hat{\tilde{w}}_t(s)=\tilde{w}(s)-\int^{s\wedge\tau}_0\nabla_v \log\alpha(t-r,w(r))dr, \ 0\leq s\leq t.
\end{equation}
Throughout the paper derivatives will be taken in $v\in G,$ so we will omit the index $v$ in the derivatives' notation.

Consider a probability measure $Q_{t,u}$ on $(\Omega,\mathcal{F}_t),$ defined via the density
$$
\frac{dQ_{t,u}}{dQ_u}=\alpha(t,u)^{-1}1_{\tau>t}.
$$
The key observation leading to the theorem \ref{thm_IW} is that on the probability space $(\Omega,\mathcal{F}_t,Q_{t,u})$ the process $\hat{\tilde{w}}_t$ is a Wiener process \cite[Ch. VIII, Th. (1.4)]{RY}.

Introduce following operators:
\begin{enumerate}
 \item $T\psi(v)=\mathbb{E}_{v} 1_{\tau<\infty} \rho(w(\tau)) \psi(w(\tau)), \ v\in G.$

Denote $\mu_v$ the distribution of $w(\tau)$ relatively to the measure $1_{\tau<\infty}d\mathbb{P}_v.$ Then the action of the operator $T$ reduces to the integration with respect to $\mu_v:$
$$
T\psi(v)=\int \psi(x)\mu_v(dx).
$$

\item $\widetilde{T}\psi(v)=\beta(v)^{-1}T\psi(v), \ v\in G.$

The operator $\widetilde{T}$ is the expectation relatively to the probability measure $Q_v:$
$$
\widetilde{T}\psi(v)=\mathbb{E}_{Q_v}\rho(w(\tau)).
$$

\item $T^k_s \psi(v)=\mathbb{E}_{Q_v} 1_{\tau>s}\psi(w(s)), \ s>0, v\in G.$

From equations \eqref{tilde_w}, \eqref{hat_tilde_w} it follows that
$$
dw(s)=\big(\nabla \log \alpha(t-s,w(s))+ \nabla \log \beta(w(s))\big)ds + d\hat{\tilde{w}}(s),
$$
where $\hat{\tilde{w}}$ is a Wiener process on $(\Omega,\mathcal{F}_t,Q_{t,u}).$ So, relatively to the measure $Q_{t,u}$ the process $w$ satisfies (degenerate) SDE. Respectively, $\{T^k_t\}_{t\geq 0}$ is the transition semigroup of a killed diffusion process $w$.
Denote $\mu_{s,v}$ the distribution of $w(s)$ relatively to the measure $1_{\tau>s}dQ_v.$ Then the action of the operator $T^k_s$ reduces to the integration with respect to $\mu_{s,v}:$
$$
T^k_s \psi(v)=\int \psi(x)\mu_{s,v}(dx).
$$

\item $\widetilde{T^k_s}\psi(v)=\alpha(s,v)^{-1}T^k_s\psi(v), \ s>0, v\in G.$

The operator $\widetilde{T^k_s}$ is the expectation relatively to the probability measure $Q_{s,v}:$
$$
\widetilde{T^k_s}\psi(v)=\mathbb{E}_{Q_{s,v}}\psi(w(s)).
$$
\end{enumerate}

The following theorem is the main result of the paper.

\begin{theorem}
\label{KV_stopped}
For every $\varphi\in L^2(\rho d \mu_u)$ the expansion \eqref{intro_IW_stopped_density} has the form
\begin{equation}
\label{KV_stopped_eq}
\begin{gathered}
\varphi(w(\tau))=\sum^\infty_{n=0}
\underset{0<t_1<\ldots<t_n}{\int\ldots\int}\alpha(t_n,u)^{-1} \bigg(\alpha(t_1,\cdot)\widetilde{T}^k_{t_1}\alpha(t_2-t_1,\cdot)\nabla \widetilde{T}^k_{t_2-t_1}\ldots \\
\alpha(t_n-t_{n-1},\cdot)\nabla \widetilde{T}^k_{t_n-t_{n-1}} \nabla \widetilde{T}\varphi\bigg)(u)
d\hat{\tilde{w}}_{t_n}(t_1)\ldots d\hat{\tilde{w}}_{t_n}(t_{n-1})d\tilde{w}(t_{n}).
\end{gathered}
\end{equation}
\end{theorem}

The proof is divided into two lemmas, which are proved in the next sections. At first we derive the Clark representation for $\varphi(w(\tau))$ with respect to the stopped Wiener process $\tilde{w}$ \cite[Ch. V, Th. (3.5)]{RY}

\begin{lemma}
\label{Clark_stopped}
For every $\varphi\in L^2(\rho d \mu_u),$ one has the representation
\begin{equation}
\label{Clark_stopped_eq}
\varphi(w(\tau))=\widetilde{T}\varphi(u)+\int^\tau_0 \nabla \widetilde{T}\varphi(w(t))d\tilde{w}(t), \ Q_u-\mbox{a.s.}
\end{equation}
\end{lemma}

Subsequently, we find the It\^o-Wiener expansion for the random variable $\psi(w(t))$ with respect to the Wiener process $\hat{\tilde{w}}.$

\begin{lemma}
\label{KV_killed}
For every $\psi\in L^2(\mu_{t,u})$ the It\^o-Wiener expansion of $\psi(w(t))$ has the form
\begin{equation}
\label{KV_killed_eq}
\begin{gathered}
\psi(w(t))=\sum^\infty_{n=0}
\underset{0<t_1<\ldots<t_n<t}{\int\ldots\int}\alpha(t,u)^{-1} \bigg(\alpha(t_1,\cdot)\widetilde{T}^k_{t_1}\alpha(t_2-t_1,\cdot)\nabla \widetilde{T}^k_{t_2-t_1}\ldots \\
\alpha(t-t_n,\cdot)\nabla \widetilde{T}^k_{t-t_n}\varphi\bigg)(u)
d\hat{\tilde{w}}_t(t_1)\ldots d\hat{\tilde{w}}_t(t_n).
\end{gathered}
\end{equation}
\end{lemma}

The theorem \ref{KV_stopped} follows by substituting $\psi=\nabla \widetilde{T}\varphi$ in \eqref{KV_killed_eq} and inserting the right-hand side of \eqref{KV_killed_eq} into \eqref{Clark_stopped_eq}.

\section{Clark Representation Formula with respect to the Measure $Q_{u}.$ Proof of the Lemma \ref{Clark_stopped}}

\begin{proof}
1) At first we will prove that the function $\widetilde{T}\varphi$ is smooth and satisfies the equation
\begin{equation}
\label{pde}
(\nabla \widetilde{T} \varphi,\nabla \log \beta)+\frac{1}{2}\Delta\widetilde{T} \varphi=0
\end{equation}
in $G.$ Indeed,
\begin{equation}
\label{T_eq}
\widetilde{T}\varphi (v)=\frac{\mathbb{E}_v 1_{\tau<\infty}\rho(w(\tau))\varphi(w(\tau))}{\beta(v)}
\end{equation}
is the ratio of two harmonic functions \cite[Ch. 4, Th. 3.7]{Port_Stone} (for the numerator the condition $\varphi\in L^2(\rho d \mu_u)$ is used). The equation \eqref{pde} is checked by straightforward calculation.

2) We will prove the relation \eqref{Clark_stopped_eq} for bounded and continuous functions $\varphi$ and $\rho$, the other cases being covered by the usual limiting procedure. Let $\{G_n\}_{n\geq 1}$ be a sequence of open relatively compact sets,  such that $\overline{G_n}\subset G$ and $G=\bigcup^\infty_{n=1} G_n.$ Denote $\tau_n$ be the  exit time from $G_n:$
$$
\tau_n=\inf\{t\geq 0: w(t)\not \in G_n\}.
$$
The convergence $\tau_n\to \tau,$ $n\to \infty,$ holds.

From the relation \eqref{tilde_w} it follows that the stopped process $w(\cdot\wedge\tau_n)$ satisfies the SDE
$$
dw(s)=\nabla\log \beta (w(s))ds+ d\tilde{w}(s), 0\leq s\leq \tau_n.
$$
Applying the It\^o formula to the function $\widetilde{T} \varphi$ and the process $w(\cdot\wedge \tau_n),$ and using \eqref{pde}, one gets the representation
$$
 \widetilde{T} \varphi(w(\tau_n))= \widetilde{T} \varphi(u)+\int^{\tau_n}_0 \nabla  \widetilde{T}\varphi (w(s))d\tilde{w}(s),  \ Q_u-\mbox{a.s.}
$$
It remains to check that $\widetilde{T} \varphi(w(\tau_n))\to \widetilde{T} \varphi(w(\tau_n)).$ As the function $\varphi$ is bounded, one has
$$
\sup_{n\geq 1} \int^\infty_0 \mathbb{E}_u1_{\tau_n>s}(\nabla  \widetilde{T}\varphi (w(s)))^2ds<\infty.
$$
Now, the convergence $\tau_n\to \tau,$ $n\to \infty,$ implies the convergence
$$
\int^{\tau_n}_0 \nabla  \widetilde{T}\varphi(w(s))d\tilde{w}(s)\xrightarrow{L^2(Q_u)}  \int^{\tau}_0 \nabla  \widetilde{T}\varphi(w(s))d\tilde{w}(s), \ n\to \infty.
$$
It remains to check that $\widetilde{T}\varphi(w(\tau_n))\to \varphi(w(\tau)),$ $n\to \infty.$ By \cite[Ch. 4, Th. 2.3]{Port_Stone} the point $w(\tau)$ is the regular point for the Dirichlet problem on $G.$
The needed convergence follows from the representation \eqref{T_eq}.
\end{proof}

\section{The Krylov-Veretennikov Formula. Proof of the Lemma \ref{KV_killed}}

\begin{proof}
The kernels $a_n$ in the expansion
$$
\psi(w(t))=\sum^\infty_{n=0}
\underset{0<t_1<\ldots<t_n<t}{\int\ldots\int}a_n(t_1,\ldots,t_n)d\hat{\tilde{w}}_t(t_1)\ldots d\hat{\tilde{w}}_t(t_n)
$$
will be recovered from the expression
$$
\begin{gathered}
\mathbb{E}_{Q_{t,u}}\psi(w(t)) \underset{0<t_1<\ldots<t_n<t}{\int\ldots\int}b_n(t_1,\ldots,t_n)d\hat{\tilde{w}}_t(t_1)\ldots d\hat{\tilde{w}}_t(t_n)=\underset{0<t_1<\ldots<t_n<t}{\int\ldots\int}\alpha(t,u)^{-1} \\
 \bigg(\alpha(t_1,\cdot)\widetilde{T}^k_{t_1}\alpha(t_2-t_1,\cdot)\nabla \widetilde{T}^k_{t_2-t_1}\ldots \alpha(t-t_n,\cdot)\nabla \widetilde{T}^k_{t-t_n}\varphi\bigg)(u)b_n(t_1,\ldots,t_n)dt_1\ldots dt_n,
\end{gathered}
$$
in which $b_n$ is a deterministic square integrable function. By induction, it is enough to check that for any square integrable $\hat{\tilde{w}}-$adapted process $\{g(s)\}_{0\leq s\leq t},$ one has
\begin{equation}
\label{eq1}
\mathbb{E}_{Q_{t,u}}\psi(w(t)) \int^t_0 g(s) d\hat{\tilde{w}}_t(s)=\int^t_0 \frac{\alpha(s,u)}{\alpha(t,u)}\mathbb{E}_{Q_{s,u}} \alpha(t-s,w(s))\nabla \widetilde{T}^k_{t-s}\psi(w(s))g(s)ds.
\end{equation}
To do it note the equalities, which follow from \eqref{hat_tilde_w} and lemma \ref{Clark_stopped}
$$
\int^t_0 g(s) d\hat{\tilde{w}}_t(s)=\int^t_0 g(s) d{\tilde{w}}(s)-\int^t_0 g(s)\nabla\log \alpha(t-s,w(s)) ds, \ Q_{t,u}-\mbox{a.s.},
$$
$$
1_{\tau>t}\psi(w(t))=T^k_t\psi(u)+\int^{t\wedge\tau}_0 \nabla T^k_{t-s}\psi(w(s))d\tilde{w}(s), \ Q_u-\mbox{a.s.}
$$
Consequently,
$$
\mathbb{E}_{Q_{t,u}}\psi(w(t)) \int^t_0 g(s) d\hat{\tilde{w}}_t(s)=
$$
$$
=\mathbb{E}_{Q_{t,u}}\psi(w(t)) \int^t_0 g(s) d{\tilde{w}}_t(s)-\mathbb{E}_{Q_{t,u}}\psi(w(t)) \int^t_0 g(s)\nabla\log \alpha(t-s,w(s)) ds=
$$
$$
=\alpha(t,u)^{-1}\bigg(\mathbb{E}_{Q_{u}}1_{\tau>t}\psi(w(t)) \int^t_0 g(s) d{\tilde{w}}_t(s)-
$$
$$
-\mathbb{E}_{Q_{u}}1_{\tau>t}\psi(w(t)) \int^t_0 g(s)\nabla\log \alpha(t-s,w(s)) ds\bigg)=
$$
$$
=\alpha(t,u)^{-1}\bigg(\int^t_0\mathbb{E}_{Q_u}1_{\tau>s}\nabla T^k_{t-s}\psi(w(s))g(s)ds-
$$
$$
-\int^t_0\mathbb{E}_{Q_u}1_{\tau>s}T^k_{t-s}\psi(w(s)) g(s)\nabla\log \alpha(t-s,w(s)) ds\bigg)=
$$
$$
=\alpha(t,u)^{-1}\int^t_0\mathbb{E}_{Q_u}1_{\tau>s}\bigg(\nabla T^k_{t-s}\psi(w(s))-T^k_{t-s}\psi(w(s))\nabla\log \alpha(t-s,w(s))\bigg) g(s)ds=
$$
$$
=\int^t_0\frac{\alpha(s,u)}{\alpha(t,u)}\mathbb{E}_{Q_{s,u}}\alpha(t-s,w(s))\nabla \widetilde{T^k_{t-s}}\psi(w(s))g(s)ds.
$$
The equality \eqref{eq1} is proved.
\end{proof}

\end{document}